\documentclass[12pt,a4paper]{amsart}
\usepackage[utf8]{inputenc}	
\usepackage[T1]{fontenc}
\usepackage[in,headings]{fullpage}
\usepackage{amsbsy, amscd, amsfonts, amsmath, amsrefs, amssymb, amsthm}
\usepackage{graphicx}
\usepackage{float}
\usepackage{color}   %May be necessary if you want to color links
\usepackage[colorlinks=true,linkcolor=blue,citecolor=blue,urlcolor=blue]{hyperref}

\newtheorem{theorem}{Theorem}%[section]

\newtheorem{corollary}[theorem]{Corollary}
\newtheorem{proposition}[theorem]{Proposition}

\theoremstyle{definition}

\theoremstyle{remark}

\newcommand{\C}{\mathbf{C}}

\renewcommand{\Re}{\mathop{\mathrm{Re}}\nolimits}
\renewcommand{\Im}{\mathop{\mathrm{Im}}\nolimits}
\newcommand{\Rzeta}{\mathop{\mathcal R }\nolimits}

\newfont{\cmbsy}{cmbsy10}
\newfont{\cmmib}{cmmib10}
\newcommand{\Orden}{\mathop{\hbox{\cmbsy O}}\nolimits}

\def\Turing{\Orden^*}

\begin{document}

\title[Simple bounds of $\Rzeta$]
{Simple bounds for the auxiliary function of Riemann.}
\author[Arias de Reyna]{J. Arias de Reyna}
\address{%
Universidad de Sevilla \\ 
Facultad de Matem\'aticas \\ 
c/Tarfia, sn \\ 
41012-Sevilla \\ 
Spain.} 

% AMS subject classifications (used in AMS journals)

% AMS keywords (used in AMS journals)

% acknowledge support, etc
%\thanks{This research was  supported by MINECO grant MTM2015--63699-P}
% \thanks{We would like to thank our colleagues for their helpful
%  criticism.}

\email{arias@us.es, ariasdereyna1947@gmail.com}

%\date{\today, \texttt{92-SimpleBounds-v5.tex}}

\begin{abstract}
We give simple numerical bounds for $\zeta(s)$, $\vartheta(s)$, $\Rzeta(s)$, $Z(t)$, for use in the numerical computation of these functions.

The purpose of the paper is to give  bounds for several functions needed in the calculation of $\zeta(s)$ and $Z(t)$.  We are not pretending to have any originality. Our object is to be useful as a  reference in our work for simple (simple to compute and simple to check) bounds. 
\end{abstract}

\maketitle

\setcounter{section}{0}
\section{Introduction.}

After publishing \cite{A1} with the aid of F. Johansson I implemented Riemann-Siegel to compute $\zeta(s)$ and related functions in mpmath. For this, we tried in \cite{Ainfty}  to perform a rigorous analysis, given $s\in\C\smallsetminus\{1\}$ and $\varepsilon>0$,  of what to compute and to what precision to obtain a complex number  $A$ such that   $|\zeta(s)-A|<\varepsilon$. For this, we need rough bound of these functions.
In fact, if we want to compute a product $ab$, and if we want the result
with an error less than $\varepsilon$, we must compute the
approximate values \texttt{a} and \texttt{b} and compute
the product $\texttt{p}=\texttt{a*b}$. Roughly  to get
$|\texttt{a*b}-ab|<\varepsilon$ we must know how big is $a$
to know the approximation to compute $b$.  On
the other hand, these bounds need not be very precise.  In fact, we want them
to determine that we must compute $b$ with more or less
$\log_2|a|$ binary digits. Hence, we prefer easy-to-compute bounds rather than tight-fitting ones.

The purpose of the paper is to give these bounds for several
functions needed in the calculation of $\zeta(s)$ and $Z(t)$.  We are not pretending to have any originality. Our object is to be useful as a  reference in our work for simple (simple to compute and simple to check) bounds. 

There are many recent works on giving explicit bounds on number theory, the lector trying to get bibliography about this may look at the web \href{http://iml.univ-mrs.fr/~ramare/TME-EMT/accueil.html}{Explicit Multiplicative Number Theory} of Olivier Ramaré. Also a good reference for inequalities  is the \href{https://fungrim.org}{Fungrim} of Fredrik Johansson. 

Here $\zeta(s)$ denotes the Riemann's zeta function. It satisfies (Titchmarsh \cite{T}*{p.~78})  the functional equation
\[\zeta(s)=\chi(s)\zeta(1-s), \quad \text{where}\quad 
\chi(s)=\frac{(2\pi)^s}{2\Gamma(s)\cos\frac{\pi s}{2}}.\]
Then $\chi(s)$ is meromorphic in $\C$ satisfying $\chi(s)\chi(1-s)=1$. 

In Titchmarsh \cite{T}*{p.~89} are introduced the Riemann-Siegel $Z$-function and the Riemann Siegel $\vartheta$-function. They are real-valued  real-analytic functions with 
\[\zeta(\tfrac12+it)=e^{-i\vartheta(t)}Z(t).\]

The function $\chi(s)$ is a meromorphic function on $\C$ with simple zeros
at the points $s=-2n$, when $n$ is a positive integer or $0$ and
poles at odd natural numbers  $s=2n+1$.  Therefore, an analytic branch of $\log\chi(s)$ can be defined in the region $\Omega$ equal to the complex numbers $\C$ with two cuts along the real axis $(-\infty,0]$ and $[1,\infty)$. We have $\chi(1/2)=1$ and we fix the branch by taking 
$\log\chi(1/2)=0$.

Then we have 
\[\vartheta(t)=-\frac{1}{2i}\log\chi(\tfrac12+it).\]
And $\vartheta(t)$ extends  analytically  on $\C$ minus two cuts along the imaginary axis $i[1/2,\infty)$ and $-i[1/2,\infty)$.

\section{Gamma function.}
The $\Gamma(s)$ function is meromorphic with simple poles at negative integers and $0$. This function never vanishes, so an analytic branch of $\log \Gamma(s)$ can be defined on the complex plane with a cut along $(-\infty,0]$. This branch is fixed by taking $\log\Gamma(s)$ real for $s>0$ real. When we write $\log\Gamma(s)$ we refer to this branch.

For $z\in\C\smallsetminus(-\infty,0]$ define 
\begin{equation}
\mu(z):=\log\Gamma(z)-\bigl(z-\tfrac12 \bigr)\log
z+z-\log\sqrt{2\pi},
\end{equation}
where both logarithms are taken as real for $z$ real and positive.

\begin{proposition} For $\Re z>0$ and $K\ge 1$ a natural number, we have 
\begin{equation}\label{asymlogGamma}
\mu(z)=\sum_{n=1}^K \frac{B_{2n}}{(2n-1)2n
z^{2n-1}}+r_K(z),
\end{equation}
where
\begin{equation}
|r_K(z)|\le \frac{|B_{2K+2}|}{(2K-1)(2K+2)|z|^{2K+1}}
\Bigl(1+\frac{2K+1}{2}\sqrt{\frac{\pi}{K}}\Bigr).
\end{equation}
\end{proposition}

We may find the proof in the book of Behnke and Sommer
\cite{BS}*{p.~304--307}.

\begin{proposition}\label{desgamma}
For $\sigma>0$ and $|s|\ge1$ we have
\begin{equation}\label{desgamaeq}
2e^{-\sigma}\; e^{-\pi
|t|/2}(\sigma^2+t^2)^{\frac{\sigma}{2}-\frac{1}{4}}
< |\Gamma(\sigma+it)|< 3\; e^{-\pi
|t|/2}(\sigma^2+t^2)^{\frac{\sigma}{2}-\frac{1}{4}}.
\end{equation}
\end{proposition}

\begin{proof} For $\sigma>0$, taking $K=1$ in
Proposition \eqref{asymlogGamma} we get
\begin{displaymath}
\log\Gamma(s)=\bigl(s-\tfrac12 \bigr)\log s-s+\log\sqrt{2\pi}
+\frac{1}{12s}+R,
\end{displaymath}
where, assuming $|s|>1$
\begin{displaymath}
|R|\le \frac{1}{30\cdot 12 |s|^3}\Bigl(1+
\frac32\sqrt{\pi}\Bigr) \le \frac{0.010164}{|s|^3}.
\end{displaymath}

Therefore,
\begin{equation}\label{logGammaapprox}
\log\Gamma(s)=\bigl(s-\tfrac12 \bigr)\log s-s+\log\sqrt{2\pi}
+\Turing\Bigl(\frac{0.0935}{|s|}\Bigr),\qquad \sigma>0,\quad
|s|>1.
\end{equation}
Taking real parts
\begin{displaymath}
\log|\Gamma(\sigma+it)|=\bigl(\sigma-\tfrac12 \bigr)\log(\sigma^2+t^2)^{1/2}
-t\arctan\frac{t}{\sigma}-\sigma+\log\sqrt{2\pi}+\Turing\Bigl(\frac{0.0935}{|s|}\Bigr).
\end{displaymath}
Assuming for the moment that $t>0$ we get
\begin{displaymath}
\log|\Gamma(\sigma+it)|
=\bigl(\sigma-\tfrac12 \bigr)\log(\sigma^2+t^2)^{1/2} -\frac{\pi
t}{2}+t\arctan\frac{\sigma}{t}-\sigma+\log\sqrt{2\pi}
+\Turing\Bigl(\frac{0.1}{|s|}\Bigr).
\end{displaymath}
For $\sigma$ and $t>0$ we have
$-\sigma<t\arctan\frac{\sigma}{t}-\sigma<0$, so that
\[
\sqrt{2\pi}\; e^{-\pi
t/2}(\sigma^2+t^2)^{\frac{\sigma}{2}-\frac{1}{4}}
e^{-\sigma-0.1/|s|} 
\le |\Gamma(\sigma+it)|\le \sqrt{2\pi}\, e^{-\pi
t/2}(\sigma^2+t^2)^{\frac{\sigma}{2}-\frac{1}{4}} e^{0.1/|s|}.
\]
(By continuity, this is also true for $t=0$.)  Since $1\le |s|$
we get $2<\sqrt{2\pi}e^{-0.1/|s|}\le \sqrt{2\pi}<3$,
so that
\begin{displaymath}
2 e^{-\pi
t/2}(\sigma^2+t^2)^{\frac{\sigma}{2}-\frac{1}{4}}e^{-\sigma} <
|\Gamma(\sigma+it)|< 3 e^{-\pi
t/2}(\sigma^2+t^2)^{\frac{\sigma}{2}-\frac{1}{4}}.
\end{displaymath}
For $t<0$ we apply that
$|\Gamma(\sigma+it)|=|\Gamma(\sigma-it)|$.
\end{proof}

\begin{proposition}\label{P:chione} 
For $\sigma>0$ and $t>1/2$ we have
\[|\chi(\sigma+it)|\le (2\pi e)^\sigma |s|^{\frac12-\sigma}.\]
For  $\sigma>0$, $t>1/2$ and $|s|\ge2\pi e$
\begin{equation}\label{chimas}
|\chi(\sigma+it)|\le (\sigma^2+t^2)^{\frac{1}{4}}.
\end{equation}
\end{proposition}

\begin{proof} We have
\begin{displaymath}
\chi(s)=\frac{(2\pi)^s}{2\Gamma(s)\cos\frac{\pi s}{2}}.
\end{displaymath}
For $t>1/2$ we have $e^{-\pi t/2}<\tfrac12  e^{\pi t/2}$ so that
\begin{displaymath}
\Bigl|\cos\frac{\pi s}{2}\Bigr|\ge \frac{e^{\pi t/2}-e^{-\pi
t/2}}{2}>\frac14 e^{\pi t/2}.
\end{displaymath}

Applying  Proposition \eqref{desgamma} for $\sigma>0$ and
$t>1/2$ we get
\[
|\chi(s)|\le\frac{ (2\pi)^{\sigma}}{4 e^{-\sigma}e^{-\pi t/2}
(\sigma^2+t^2)^{\frac{\sigma}{2}-\frac{1}{4}}
\frac{1}{4}e^{\pi t/2}}= (2\pi
e)^{\sigma}(\sigma^2+t^2)^{\frac{1}{4}-\frac{\sigma}{2}}
= (\sigma^2+t^2)^{\frac{1}{4}}\Bigl(\frac{4\pi^2
e^2}{\sigma^2+t^2}\Bigr)^{\frac{\sigma}{2}}.
\]
Therefore, we have
\begin{equation}\label{chimenos2}
|\chi(s)|\le (\sigma^2+t^2)^{\frac{1}{4}}\Bigl(\frac{4\pi^2
e^2}{\sigma^2+t^2}\Bigr)^{\frac{\sigma}{2}}, \qquad (\sigma>0,\ t>\tfrac12 ).
\end{equation}
For $|s|\ge2\pi e$  we get \eqref{chimas}.
\end{proof}

\begin{proposition} For $\sigma\le0$ and $t\ge\tfrac12 $ we have
\begin{equation}\label{chimenos}
|\chi(\sigma+it)|\le
\frac{6}{(2\pi)^{1-\sigma}}\{(1-\sigma)^2+t^2\}^{\frac14-\frac{\sigma}{2}}.
\end{equation}
\end{proposition}

\begin{proof}
We have $\chi(s)\chi(1-s)=1$ for all $s$.  Hence, we have
\[
|\chi(s)|=\frac{1}{|\chi(1-s)|} = \Bigl|\frac{2\Gamma(1-s)\cos\frac{\pi(1-s)}{2}}{(2\pi)^{1-s}}\Bigr|=\Bigl|\frac{2\Gamma(1-s)\sin\frac{\pi s}{2}}{(2\pi)^{1-s}}\Bigr|.
\]
For $\sigma\le0$ we have $|1-s|\ge1$ and by \eqref{desgamaeq} 
\[|\Gamma(1-s)|\le 3 e^{-\pi|t|/2}((1-\sigma)^2+t^2)^{\frac14-\frac{\sigma}{2}},\]
and 
\[\Bigl|\sin\frac{\pi s}{2}\Bigr|\le e^{\pi|t|/2}.\]
Therefore,
\[|\chi(s)|\le \frac{6}{(2\pi)^{1-\sigma}}((1-\sigma)^2+t^2)^{\frac14-\frac{\sigma}{2}}.\qedhere\]
\end{proof}

\section{Bounds for \texorpdfstring{$\vartheta(t)$}{Riemann Siegel theta}.}

\begin{proposition}
We have for complex $t$
\begin{equation}
|\vartheta(t)| \le 2|t|\log|t|,\qquad |t|\ge4,\quad |\Re t|\ge 1.
\end{equation}
\end{proposition}

\begin{proof}
As an analytic function 
$\vartheta(t)=-\frac{1}{2i}\log\chi(\frac12+it)$ extends  analytically  on $\C$ minus two cuts along the imaginary axis $i[1/2,\infty)$ and $-i[1/2,\infty)$.

Hence, for $t$ complex and $s=\frac12+it$ we have
\begin{multline*}
\vartheta(t)=-\frac{1}{2i}\log\chi(s)=-\frac{1}{2i}\log\frac{(2\pi)^s}{2\cos\frac{\pi
s}{2}\,\Gamma(s)}=\\= -\frac{1}{2i}\Bigl\{s\log2\pi-\log\bigl(
e^{-\pi i s/2}(1+e^{\pi i s})\bigr)-\log\Gamma(s)\Bigr\},
\end{multline*}
for an adequate election of the logarithms. 

For $x:=\Re s>0$ and $y:= \Im s\ge 1$  we have  $|e^{\pi i s}|\le e^{-\pi y}\le e^{-\pi}$, therefore
\[\log\bigl(
e^{-\pi i s/2}(1+e^{\pi i s})\bigr)=-\frac{\pi i s}{2}+\Turing\Bigl(\frac{e^{-\pi y}}{1-e^{-\pi}}\Bigr)=-\frac{\pi i s}{2}+\Turing\Bigl(\frac{0.05}{y}\Bigr).\]
Since we assume that $x=\Re s>0$ and $|s|\ge y\ge1$ it follows by \eqref{logGammaapprox} that
\begin{align*}
\vartheta(t)&= -\frac{1}{2i}\Bigl\{s\log2\pi+\frac{\pi
is}{2}+\Turing\bigl(\tfrac{0.05}{y}\bigr)-\bigl(s-\tfrac12 \bigr)\log
s+s-\log\sqrt{2\pi}
+\Turing\bigl(\tfrac{0.0935}{|s|}\bigr)\Bigr\}\\
&=\frac{s-\frac12}{2i}\log\frac{s}{2\pi}-\frac{\pi s}{4}+\frac{i s}{2}+\Turing\bigl(\tfrac{0.15}{y}\bigr).
\end{align*}
It is easy to check that the elections taken of the logarithms make the right-hand side here almost real for real $t\gg1$, therefore there is no multiple of $2\pi i$ to add in this equality.

Therefore, for  $s=\frac12+it$, $s=x+iy$ with $x>0$, $y>1$ and $|t|>1$
\begin{multline*}
|\vartheta(t)|\le
\frac{|t|}{2}\log\frac{|t|+1/2}{2\pi}+\frac{\pi
|t|}{4}+\frac{\pi
|t|}{4}+\frac{\pi}{8}+\frac{|t|}{2}+\frac{1}{4}+0.15\le\\  \le
\frac{|t|}{2}\log\frac{3|t|}{4\pi}+\frac{\pi
|t|}{2}+\frac{\pi}{8}+\frac{|t|}{2}+0.4\le
\frac{|t|}{2}\log|t|+1.36\,|t|+0.80.
\end{multline*}
And for $|t|>4$ we will have
\begin{equation}\label{boundvartheta}
|\vartheta(t)| \le 2|t|\log|t|
\end{equation}
Since $x=\Re s=\tfrac12 -\Im t$ and $y=\Im s = \Re t$, the
inequality \eqref{boundvartheta} is true for $|t|>4$, $\Re t>1$
and $\Im t <\tfrac12 $. Since $\vartheta(t)$ is real for real $t$, it takes complex conjugate values in conjugate arguments. So the inequality is also true for $\Im t>-\tfrac12 $.
It follows that the inequality is true for $|t|>4$ and $\Re t>1$.  Since $\vartheta(-t)=-\vartheta(t)$, 
the inequality is also true for $|t|>4$, $|\Re t|>1$.
\end{proof}

\section{Bounds for \texorpdfstring{$\zeta(s)$}{zeta}.}
In this section assume that $s=\sigma+it$ where $\sigma =\Re s$ and
$t=\Im s$.
\begin{proposition}
We have
\begin{equation}\label{zetamas2}
|\zeta(\sigma+it)|<2,\qquad \sigma\ge2.
\end{equation}
\end{proposition}

\begin{proof}
For $\sigma\ge 2$ we have $|\zeta(\sigma+it)|\le
|\zeta(\sigma)|\le \zeta(2)< 2.$
\end{proof}

\begin{proposition}
We have
\begin{equation}
|\zeta(\sigma+it)|\le 2(2\pi)^{\sigma}\{(1-\sigma)^2+t^2\}^{\frac14-\frac\sigma2}
,\qquad \sigma\le -1,\  |t|\ge\tfrac12 .
\end{equation}
\end{proposition}

\begin{proof}
We apply the functional equation $\zeta(s)=\chi(s)\zeta(1-s)$.
For $\sigma\le-1$ and assuming $t\ge\tfrac12 $ by \eqref{chimenos} and
\eqref{zetamas2} we get
\begin{displaymath}
|\zeta(\sigma+it)|\le 2
\frac{6}{(2\pi)^{1-\sigma}}\{(1-\sigma)^2+t^2\}^{\frac{1}{4}-\frac{\sigma}{2}},
\end{displaymath}
and note that $6<2\pi$.
\end{proof}

\begin{proposition}\label{zetacritica}
We have
\begin{equation}
|\zeta(s)|\le 1+\frac{t}{\sigma},\qquad 0<\sigma\le2,\quad t\ge2
\end{equation}
\end{proposition}

\begin{proof} We start with  formula  (2.1.4) of Titchmarsh
\cite{T}
\begin{equation}
\zeta(s)=s\int_1^{+\infty}\frac{\lfloor x\rfloor
-x+\frac{1}{2}}{x^{s+1}}\,dx+\frac{1}{s-1}+\frac{1}{2},\qquad \sigma>0.
\end{equation}
Then,  we find for $t\ge2$
\[
|\zeta(s)|\le
|s|\int_1^{+\infty}\frac{\frac{1}{2}}{x^{\sigma+1}}\,dx+\frac{1}{|s-1|}+\frac{1}{2}\le (t+2)\frac{1}{2\sigma}+\frac{1}{t}+\frac{1}{2}\le
\frac{t}{\sigma}+1.\qedhere
\]
\end{proof}

\begin{corollary} We have
\begin{equation}\label{zetaright}
|\zeta(\sigma+it)|\le 3t,\qquad \sigma\ge\tfrac12 ,\quad t\ge2.
\end{equation}
\end{corollary}

\begin{proof} By  Proposition \ref{zetacritica} and  hypothesis $\sigma\ge1/2$
and $t\ge2$  we have (assuming first $\sigma\le2$)
\begin{displaymath}
|\zeta(\sigma+it)|\le1+\frac{t}{\sigma}\le 1+2t\le 3t.
\end{displaymath}
For $\sigma\ge2$ the result follows from \eqref{zetamas2}.
\end{proof}

\section{Function \texorpdfstring{$\Rzeta(s)$}{R(s)}.}

\begin{theorem}\label{rzetazetasum}
We have
\begin{equation}
\Bigl|\Rzeta(s)-\sum_{n\le \sqrt{t/2
\pi}}\frac{1}{n^s}\Bigr|\le
\Bigl(\frac{t}{2\pi}\Bigr)^{-\sigma/2}, \qquad 0\le\sigma\le1, \quad t\ge 3\pi,
\end{equation}
and
\begin{equation}
\Bigl|\Rzeta(s)-\sum_{n\le \sqrt{t/2
\pi}}\frac{1}{n^s}\Bigr|\le
\Bigl(\frac{t}{2\pi}\Bigr)^{-1/2}, \qquad \sigma\ge1, \quad t\ge 16\pi.
\end{equation}
\end{theorem}

\begin{proof}
In \cite{A1}*{Theorem 3.1} we find  the Riemann-Siegel
expansion of the integral, giving $\Rzeta(s)$. There it is
proved that for $t>0$ and $\sigma>0$
\begin{equation}
\Rzeta(s)=\sum_{n=1}^N\frac{1}{n^s}+(-1)^{N-1}U
a^{-\sigma}\{ C_0(p)+\frac{C_1(p)}{a}+RS_1\},
\end{equation}
where $a=\sqrt{t/2\pi}$, $N=\lfloor a\rfloor$, $p=1-2(a-N)$ and
\begin{equation}
U=\exp\Bigl\{-i\Bigl(\frac{t}{2}\log\frac{t}{2\pi}-\frac{t}{2}-
\frac{\pi}{8}\Bigr)\Bigr\},
\end{equation}
also  \cite{A1}*{(5.2)}
\begin{equation}
C_0(p)=F(p):=\frac{e^{\pi i
(\frac{p^2}{2}+\frac38)}-i\sqrt{2}\cos\frac{\pi}{2}p}{2\cos\pi
p}.
\end{equation}
By \cite{A1}*{(2.4)} we have the following
\begin{equation}
C_1(p)=\frac{1}{\pi^2}(d^{(1)}_0 F'''(p)+ d^{(1)}_1 F'(p)).
\end{equation}
By the recurrence equation for the coefficients $d^{(k)}_j$
\cite{A1}*{(2.7), (2.8), (2.9) and (2.10)}
\begin{equation}
C_1(p)=\frac{1}{\pi^2}(\frac{1}{12} F'''(p)-\frac{\pi}{2i}(\sigma-\tfrac12 ) F'(p)).
\end{equation}
So, by \cite{A1}*{(6.7)}
\begin{equation}
|C_1(p)|\le \frac{1}{\pi^2}(\frac{1}{12} 2\pi+\frac{\pi}{2}|\sigma-\tfrac12 |)=
\frac{1}{6\pi}+\frac{|\sigma-\frac12|}{2\pi}.
\end{equation}
Also  by \cite{A1}*{Theorem 4.2}
\begin{equation}
|RS_1|\le
\frac17 2^{3\sigma/2}\frac{1}{(a/1.1)^2},
\end{equation}
and in \cite{A1}*{Theorem 6.1} it is proved  that $C_0(p)=F(p)$
satisfies $|F(p)|\le \tfrac12 $.

Therefore, the difference $D:=\Rzeta(s)-\sum_{n=1}^N n^{-s}$  in
absolute value is bounded by
\begin{equation}\label{firstinequ}
|D|\le a^{-\sigma}\Bigl\{\frac12+\frac{1}{6\pi
a}+\frac{|\sigma-\frac12|}{2\pi a} + \frac17
2^{3\sigma/2}\frac{1}{(a/1.1)^2}\Bigr\}.
\end{equation}
For $0\le \sigma\le 1$ we get
\begin{displaymath}
|D|\le a^{-\sigma}\Bigl\{\frac12+\frac{1}{6\pi
a}+\frac{1}{4\pi a} + \frac17
2^{3/2}\frac{1}{(a/1.1)^2}\Bigr\}.
\end{displaymath}
This is less than $a^{-\sigma}$ for $a\ge\sqrt{3/2}$, hence for
$t\ge 3\pi$.

For $\sigma\ge1$  by \eqref{firstinequ} we get
\begin{displaymath}
|D|\le \frac{1}{2a}+\frac{1}{6\pi
a^2}+\frac{|\sigma-\frac12|}{2\pi a \cdot a^\sigma} + \frac17
(2^{3/2}/a)^\sigma\frac{1}{(a/1.1)^2}.
\end{displaymath}
When $t\ge16\pi$, we get $2^{3/2}/a \le 1$ and
$2\sigma<2^{3\sigma/2}\le a^\sigma$
\begin{displaymath}
|D|\le \frac{1}{a}\Bigl\{\frac{1}{2}+\frac{1}{6\pi
a}+\frac{1}{2\pi} + \frac{1.21}{7a}\Bigr\}.
\end{displaymath}
and this is less than $1/a$ when $a>1$.
\end{proof}

\begin{proposition} For $\sigma\ge 2$ and $t>16\pi$ we have
\begin{equation}
|\Rzeta(s)-1|\le
\frac{3}{2^\sigma}+\Bigl(\frac{2\pi}{t}\Bigr)^{\min(\sigma,1)/2}.
\end{equation}
\end{proposition}

\begin{proof} By Theorem  \ref{rzetazetasum} we have
\begin{align*}
|\Rzeta(s)-1|&\le \sum_{n=2}^\infty
\frac{1}{n^\sigma}+\Bigl(\frac{2\pi}{t}\Bigr)^{\min(\sigma,1)/2}
\le \frac{1}{2^\sigma}+\frac{1}{3^\sigma}+
\int_3^{+\infty}\frac{dv}{v^\sigma}+\Bigl(\frac{2\pi}{t}\Bigr)^{\min(\sigma,1)/2}\\
&\le \frac{1}{2^\sigma}+\frac{4}{3^\sigma}
+\Bigl(\frac{2\pi}{t}\Bigr)^{\min(\sigma,1)/2}<\frac{3}{2^\sigma}+\Bigl(\frac{2\pi}{t}\Bigr)^{\min(\sigma,1)/2}\quad \text{(for $\sigma\ge2$)}.\qedhere
\end{align*}
\end{proof}

\begin{proposition} We have 
\begin{equation}\label{Rright}
|\Rzeta(\sigma+it)|\le 2\sqrt{\frac{t}{2\pi}}\qquad \sigma>0, \quad t>16\pi.
\end{equation}
\end{proposition}

\begin{proof}
We apply Theorem \ref{rzetazetasum}. Since we assume that
$\sigma>0$ we have $|n^{-s}|\le 1$. So the sum is bounded by
the number of terms $\sqrt{t/2\pi}$.

We also have
\begin{displaymath}
\Bigl(\frac{t}{2\pi}\Bigr)^{-\min(\sigma,1)/2}\le1\le
\Bigl(\frac{t}{2\pi}\Bigr)^{1/2} \end{displaymath} since
$t>2\pi$. This proves the result.
\end{proof}

\begin{proposition}
 For $\sigma\le 0$ and $t\ge 16\pi$ we have
\begin{equation}\label{Rzetaleft}
|\Rzeta(\sigma+it)|\le
\frac{19\,t}{(2\pi)^{1-\sigma}}\{(1-\sigma)^2+t^2\}^{\frac{1}{4}-\frac{\sigma}{2}},\qquad
\sigma\le 0,\quad t\ge 16\pi.
\end{equation}
\end{proposition}

\begin{proof} In \cite{A2} it is shown that 
\begin{displaymath}
\zeta(s)=\Rzeta(s)+\chi(s)\overline{\Rzeta}(1-s)=\Rzeta(s)+\chi(s)\overline{\Rzeta(1-\overline{s})}
\end{displaymath}
Putting $1-\sigma+it$ instead of $s$
\begin{displaymath}\zeta(1-\sigma+it)=\Rzeta(1-\sigma+it)
+\chi(1-\sigma+it)\overline{\Rzeta(\sigma+it)}
\end{displaymath}
Since $\chi(s)\chi(1-s)=1$ it follows that
\begin{displaymath}
|\Rzeta(\sigma+it)|=|\chi(\sigma-it)\{\zeta(1-\sigma+it)-\Rzeta(1-\sigma+it)\}|
\end{displaymath}
We know that $\chi(\sigma-it)=\overline{\chi(\sigma+it)}$.
Therefore,  by \eqref{chimenos}, \eqref{zetaright} and \eqref{Rright}
we find that  for $\sigma<0$ and $t>16\pi$
\begin{align*}
|\Rzeta(\sigma+it)|&\le
\frac{6}{(2\pi)^{1-\sigma}}\{(1-\sigma)^2+t^2\}^{\frac{1}{4}-\frac{\sigma}{2}}\Bigl(3t+
2\sqrt{\frac{t}{2\pi}}\Bigr)\\
&\le\frac{19t}{(2\pi)^{1-\sigma}}\{(1-\sigma)^2+t^2\}^{\frac{1}{4}-\frac{\sigma}{2}}.\qedhere
\end{align*}
\end{proof}

\begin{proposition}
For $\sigma\ge0$, $|s-1|\ge2$ and $0<x\le |s|$ we have 
\begin{equation}
\Bigl|\zeta(s)-\sum_{n\le x}\frac{1}{n^s}\Bigr|\le 7|s|x^{-\sigma}.
\end{equation}

For $\sigma>1$ we have the alternative bound
\begin{equation}\label{E:boundzeta2}
\Bigl|\zeta(s)-\sum_{n\le x}\frac{1}{n^{s}}\Bigr|\le
\frac{x^{1-\sigma}}{\sigma-1}\Bigl(1+\frac{\sigma-1}{x}\Bigr).
\end{equation}

\end{proposition}
\begin{proof}
Assume first that $\sigma\ge1/2$, then 
\[\zeta(s)=\sum_{n\le x}\frac{1}{n^s}+\frac{x^{1-s}}{s-1}+\frac{\{x\}-\frac12}{x^s}+s\int_x^{+\infty}\frac{\frac12-\{u\}}{u^{s+1}}\,du.\] 
Hence, 
\[\Bigl|\zeta(s)-\sum_{n\le x}\frac{1}{n^s}\Bigr|\le \frac{x^{1-\sigma}}{|s-1|}+\frac{1}{2}x^{-\sigma}+\frac{|s|}{2}\int_x^\infty u^{-\sigma-1}\,du\le \Bigl|\frac{s}{s-1}\Bigr|x^{-\sigma}+\frac12 x^{-\sigma}+\frac{|s|}{2\sigma}x^{-\sigma}.\]
Since $|s-1|\ge2$ we have $|s|\ge1$ and
\[\Bigl|\zeta(s)-\sum_{n\le x}\frac{1}{n^s}\Bigr|\le \Bigl(\frac12+\frac12+1\Bigr)|s|x^{-\sigma}\le 2|s|x^{-\sigma}.\]
For $0\le\sigma\le \frac12$ and $t>16\pi$ we have by \eqref{Rright} and Proposition  \ref{P:chione} 
\begin{align*}
|\zeta(s)|&\le |\Rzeta(s)+\chi(s)\overline{\Rzeta(1-\overline{s})}|\le2\sqrt{\frac{t}{2\pi}}+2\sqrt{2\pi e}|s|^{\frac12-\sigma}\sqrt{\frac{t}{2\pi}}\\
&\le \frac{2}{\sqrt{2\pi}}|s|^{1/2}+2\sqrt{e}|s|^{\frac12-\sigma}|s|^{1/2}
\le \Bigl(\frac{2}{\sqrt{2\pi}}+2\sqrt{e}\Bigr)|s|x^{-\sigma}\le 5|s|x^{-\sigma}.
\end{align*}
And a numerical study shows that $|\zeta(s)|\le 2|s|^{1/2}$ for $s\in[0,1/2]\times[\sqrt{3},16\pi]$. Assuming that $t>0$ and $|s-1|\ge2$  we have $t>\sqrt{3}$. Therefore $|\zeta(s)|\le 5|s|x^{-\sigma}$ for $|s-1|\ge2$ and $0<\sigma\le 1/2$. 

We also have for $0\le\sigma\le 1/2$ and $0<x\le |s|$
\[\Bigl|\sum_{n\le x}\frac{1}{n^s}\Bigr|\le 1+\int_1^xu^{-\sigma}\,du\le 1+\frac{x^{1-\sigma}-1}{1-\sigma}\le \frac{x^{1-\sigma}}{1-\sigma}\le 2|s| x^{-\sigma}.\]
It follows that for $0\le \sigma\le1/2$ and $|s-1|\ge2$ we have 
\[\Bigl|\zeta(s)-\sum_{n\le x}\frac{1}{n^s}\Bigr|\le7|s|x^{-\sigma}.\]

In the case $\sigma>1$ the value of $\zeta(s)$ is given by the series, so  if $m\le x<m+1$ we have 
\begin{align*}
\Bigl|\zeta(s)-\sum_{n\le x}\frac{1}{n^{s}}\Bigr|&=\Bigl|\sum_{n>x}\frac{1}{n^s}\Bigr|
\le \frac{1}{(m+1)^\sigma}+\int_{m+1}^\infty u^{-\sigma}\,du\\
&=\frac{1}{(m+1)^\sigma}+\frac{(m+1)^{1-\sigma}}{\sigma-1}<\frac{x^{1-\sigma}}{\sigma-1}\Bigl(1+\frac{\sigma-1}{x}\Bigr).\qedhere
\end{align*}

\end{proof}

\end{document}